\documentclass{amsart}
\usepackage{amssymb,amsmath,latexsym,amscd,amsfonts}
\usepackage{epsfig}
\usepackage{graphicx}
\usepackage{graphics}
\usepackage{psfrag}

\newtheorem{thm}{Theorem}

\newtheorem{cor}{Corollary}

\begin{document}

\title{A Million-Dollar Proof\footnote{Published in \emph{Mathematical Intelligencer,} vol.~29 no.~4 (2007), pg.~8.}}
\author{Aaron Abrams}
\maketitle

I came up with the main theorem of this paper during the
summer of 1992 while I was an undergraduate student at Joe Gallian's
REU at the University of Minnesota, Duluth.  Until recently, I neither
shared it with anyone nor fully appreciated its value.  But now Fermat's
Last Theorem has a proof, and it seems like the Poincar\'e Conjecture
does too, and with ever more attention directed toward problems like
these, I thought I'd better get this in print.

Here is the theorem:

\begin{thm} $N=1$.
\end{thm}

The proof is based on what I thought at the time was a surprising observation,
namely the following identity:

\begin{equation}\label{obs}
\sum\limits_{N=0}^{\infty} \frac 1 {2^N} = \sum\limits_{N=0}^{\infty} \frac N {2^N}.
\end{equation}

Now, this is almost surely false, since for one thing the terms on the right
are larger than the corresponding terms on the left.  Nonetheless there are 
many proofs of this result, most involving evaluating both sides of the equation.
(In most proofs I've seen, the common value is 2.)  In fact, one need not
evaluate either side to prove the result; moreover, one needs no words:

\begin{figure}[h]
\psfrag{a}{$1$}
\psfrag{b}{$1/2$}
\psfrag{c}{$1/4$}
\psfrag{d}{$1/8$}
\psfrag{e}{$1/16$}
\psfrag{f}{$2/4$}
\psfrag{g}{$3/8$}
\psfrag{h}{$4/16$}
\psfrag{i}{$\frac 5 {32}$}
\psfrag{j}{$2$}
\psfrag{k}{$3$}
\psfrag{l}{$4$}
\psfrag{m}{$5$}
\psfrag{n}{$6$}
\scalebox{.49}{
	\includegraphics{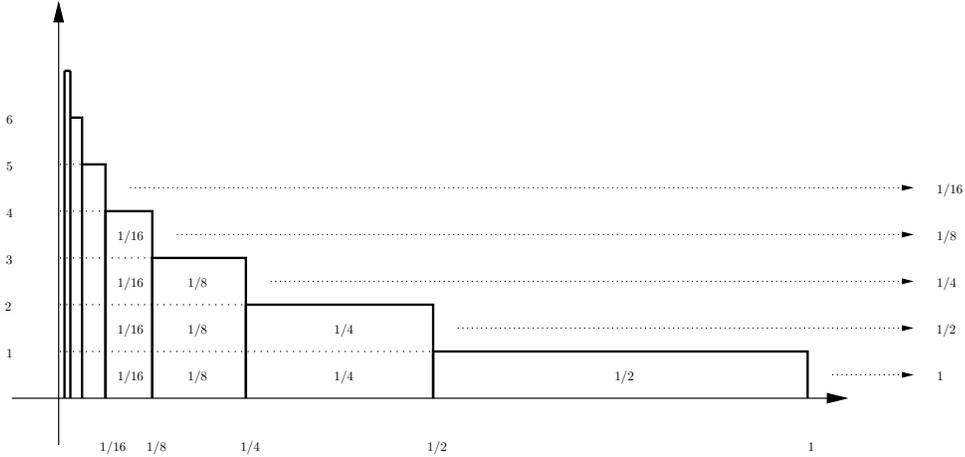}
}
\caption{$1 + \frac 1 2 + \frac 1 4 + \frac 1 8 + \frac 1 {16} + \cdots = 
\frac 1 2 + \frac 2 4 + \frac 3 8 + \frac 4 {16} + \frac 5 {32} + \cdots$.
\newline 
(Note that the $x$-axis scale is different from the $y$-axis scale.)}
\end{figure}

Anyway, beginning with (\ref{obs}), I was naturally tempted to cancel 
the sigmas and 
the denominators and thereby conclude that $N=1$.  This didn't seem quite 
rigorous, however.  Luckily Eric Wepsic was around to help me formalize 
the argument.  Here's the proof:

\begin{proof}
Define the function $f(x)=\sum\limits_{N=0}^{\infty} \frac x {2^N}$.  
Note the crucial fact that $f$ is one-to-one.

By (\ref{obs}), $f(1)=f(N)$.  Thus $N=1$.
\end{proof}

I couldn't help but notice the corollary:

\begin{cor}
$P=NP$.
\end{cor}

\end{document}